\documentclass[12pt]{article}
\newcommand{\copyleft}{
GNU FDL\thanks{
Copyright (C) 1988 Peter G. Doyle.
Permission is granted to copy, distribute and/or modify this document
under the terms of the GNU Free Documentation License, 
as published by the Free Software Foundation;
with no Invariant Sections, no Front-Cover Texts, and no Back-Cover Texts.
}}
\title{Electric currents in infinite networks}
\author{Peter G. Doyle}
\date{Version dated 25 October 1988
\\ \copyleft
}

\usepackage{epsfig}
\newcommand{\fig}[2]{
\begin{figure}
\psfig{figure=figures/#1.ps,width=370pt}
\caption{#2}
\label{#1}
\end{figure}
}
\newcommand{\energy}{{\rm energy}}
\newcommand{\resistance}{{\rm resistance}}
\newcommand{\oddresistance}{{\rm odd\_resistance}}
\newcommand{\evenresistance}{{\rm even\_resistance}}
\newcommand{\shortnetwork}{{\rm short\_network}}
\newcommand{\cutnetwork}{{\rm cut\_network}}
\newcommand{\shortgraph}{{\rm short\_graph}}
\newcommand{\cutgraph}{{\rm cut\_graph}}

\newcommand{\averageevenresistance}{{\rm average\_even\_resistance}}
\newcommand{\averageoddresistance}{{\rm average\_odd\_resistance}}
\newcommand{\averagevalence}{{\rm average\_valence}}
\newcommand{\boundary}{{\rm boundary}}
\newcommand{\edgeboundary}{{\rm edge\_boundary}}
\newcommand{\valence}{{\rm valence}}
\newcommand{\valenceone}{{\rm valence1}}
\newcommand{\valencetwo}{{\rm valence2}}

\newcommand{\og}{{O_{G}}}
\newcommand{\ohp}{{O_{HP}}}
\newcommand{\ohb}{{O_{HB}}}
\newcommand{\ohd}{{O_{HD}}}

\newcommand{\qed}{\heartsuit}
\newcommand{\dispqed}{\;\;\qed}

\begin{document}

\maketitle

\section{Introduction.}

In this survey, we will present the basic facts about conduction
in infinite networks.
This survey is based on the work of
Flanders
\cite{flanders:infinite,flanders:infiniteGrid},
Zemanian
\cite{zemanian:infinite},
and Thomassen
\cite{thomassen:infinite},
who developed the theory of infinite networks from scratch.
Here we will show how to get a more complete theory
by paralleling the well-developed theory of conduction
on open Riemann surfaces.
Like Flanders and Thomassen,
we will take as a test case for the theory
the problem of determining the resistance
across an edge of a $d$-dimensional grid of $1$ ohm resistors.
(See Figure 1.)
We will use our borrowed network theory
to unify, clarify and extend their work.

\fig{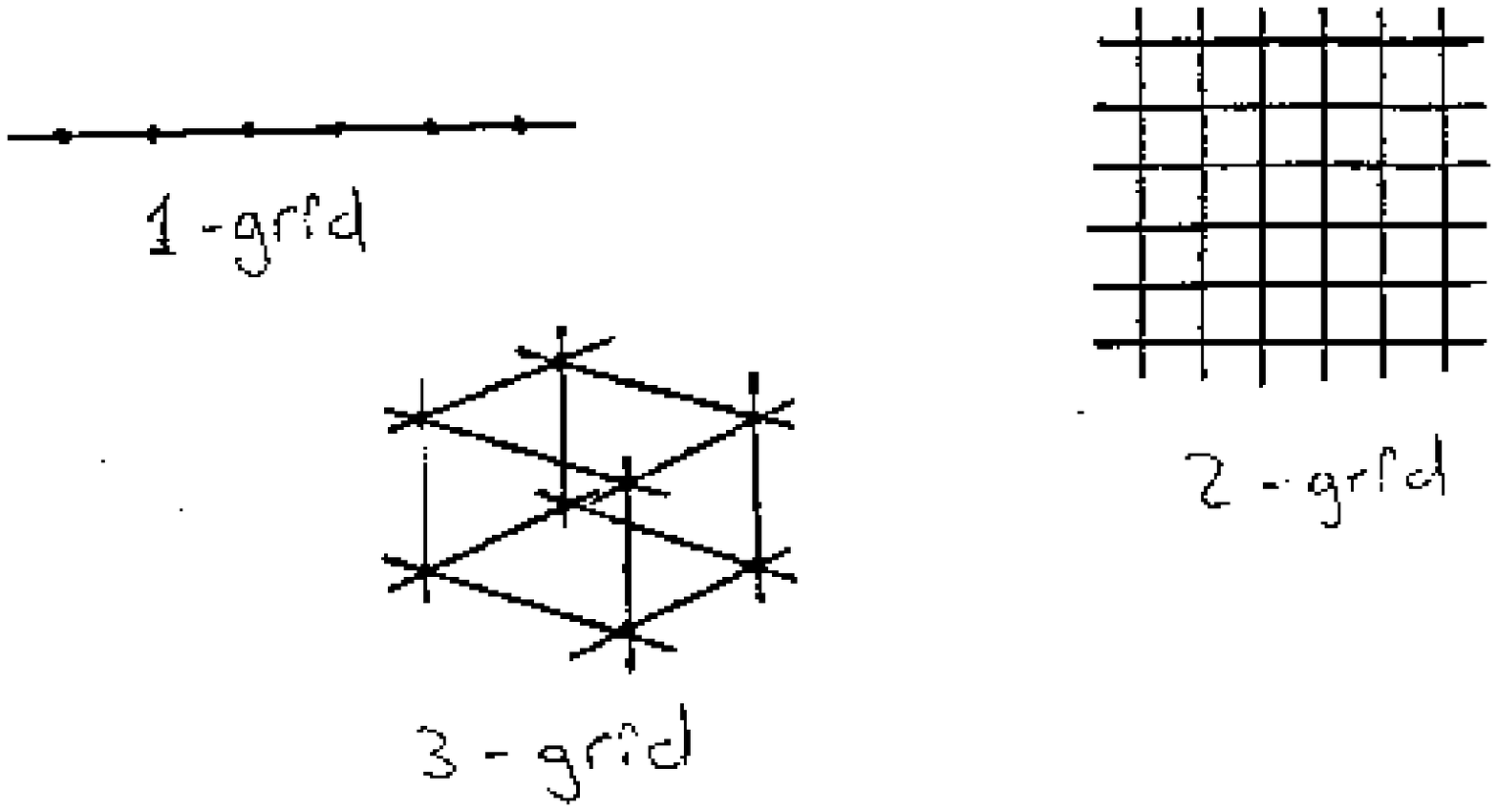}{Grids.}

\section{The engineers and the grid.}

Engineers have long known how to compute the resistance across
an edge of a $d$-dimensional grid of $1$ ohm resistors
using only the principles of symmetry and superposition:
Given two adjacent vertices $p$ and $q$,
the resistance across the edge from $p$ to $q$
is the voltage drop along the edge
when a $1$ amp current is injected at $p$ and withdrawn at $q$.
Whether or not $p$ and $q$ are adjacent,
the unit current flow from $p$ to $q$
can be written as the superposition
of the unit current flow from $p$ out to infinity
and the unit current flow from infinity into $q$.
By symmetry,
in the unit current flow from $p$ to infinity,
the flow out of $p$ is distributed equally
among the $2d$ edges going out of $p$,
so the flow along any one of them is $1 \over 2d$ amps.
Similarly, in the unit current flow into $q$,
the flow along each edge coming into $q$ is $1 \over 2d$.
If $p$ and $q$ are adjacent,
when the two flows are superimposed the flow along the edge from
$p$ to $q$ is $2 {{1}\over{2d}} = {{1}\over{d}}$
But since this edge has resistance $1$ ohm,
the voltage drop along it is also $1/d$,
so the effective resistance between $p$ and $q$ is $1/d$.

A second justification for the answer $1/d$
was offered by Foster
\cite{foster:average},
based on his theorem that
the average of the resistances across all the edges of
a finite graph is
\[
{{n-1} \over {e}}
,
\]
where $n$ is the number of vertices and $e$ the number of edges 
of the graph.
(See Foster
\cite{foster:average,foster:extension},
Weinberg
\cite{weinberg:network} pp. 170--176.)
This theorem is an immediate consequence of Kirchhoff's Rule,
also known as Maxwell's Rule,
according to which
the resistance across an edge of a finite connected
network of $1$ ohm resistors
is the probability that a randomly selected spanning tree of the
network contains that edge.
(See Kirchhoff
\cite{kirchhoff:gleichungen},
Maxwell
\cite{maxwell:treatise} pp. 409--410,
Bollob\'{a}s
\cite{bollobas:graph},
Thomassen
\cite{thomassen:infinite}.)
If we confidently apply Foster's theorem to the $d$-grid,
we get the expected answer $1/d$ for the resistance across an edge.

These two arguments could hardly be more elegant.
Mathematically speaking, though, they leave something to be desired.
As Flanders
\cite{flanders:infinite}
points out,
for an infinite network it is not even clear
what it means to talk about the effective resistance
between two nodes.
For a finite network,
Kirchhoff's laws determine a unique unit current
flow from $p$ to $q$,
but this clearly isn't true for an infinite network.
The engineers believe that in an infinite grid
there will still be a uniquely determined unit current flow
from $p$ to $q$,
and that this flow can be gotten by superimposing uniquely determined
symmetrical flows from $p$ to infinity and from infinity to $q$.
They offer various arguments to justify their beliefs,
based on approximating the infinite network by finite networks,
and the like.
The arguments I have heard are not rigorous,
and I am convinced that the true reason for the engineers' beliefs
is the analogy between conduction in the $d$-grid and classical
potential theory in $d$-dimensional Euclidean space.
Any property of conduction in Euclidean space should continue to hold
for the grid,
unless there is some obvious reason why it should not.

The engineers are right, of course.
There is only one sensible definition
of the resistance across an edge of the $d$-grid,
or rather,
there are two sensible definitions,
which turn out to agree.
The answer $1/d$ is correct,
and the arguments they give to justify it are basically sound.
This correctness of the answer $1/d$ was shown
by Flanders
\cite{flanders:infiniteGrid}
in the case $d=2$,
using the symmetry-and-superposition method;
Thomassen
\cite{thomassen:infinite}
proved the result in general,
using Foster's method.
Once we have developed the theory of infinite networks,
we will go over the work of Flanders and Thomassen,
and show just how right the engineers are.

\section{Networks as Riemann surfaces.}

We can think of an infinite network, such as the $d$-grid,
as a degenerate Riemann surface,
where all of the tubes are very long and skinny.
One indication of this is that when
asked to consider conduction in the
three-dimensional grid, many people think immediately of
a jungle gym, even though it is made of hollow tubes.
Thus, if we need to develop the theory of infinite networks,
we may hope to get the theory as a limiting case of
the theory of conduction on Riemann surfaces,
or at least to be able to develop the network theory
along the same lines as the Riemann surface theory.
It is the latter approach that we will follow here.
The network theory we will develop follows the Riemann surface
theory as described in Ahlfors and Sario's standard text
\cite{ahlforsSario:riemannSurfaces};
see also Rodin and Sario
\cite{rodinSario:principal}.
We will not duplicate any of the proofs.
For network terminology and the basic theory of finite networks
see Doyle and Snell
\cite{doyleSnell:carus},
or the beautiful and concise treatment by Thomassen
\cite{thomassen:infinite}.

\section{Infinite networks.}

Consider a connected infinite electric network,
that is,
a connected infinite graph where every edge is assigned
a resistance.
We will assume that the network is {\em locally finite},
that is,
that each vertex has finite valence.
Multiple edges connecting two distinct vertices are allowed.
Any graph can be treated as
a network by assigning each of its edges a resistance of $1$ ohm.

We will often want to approximate an infinite network
by a sequence of larger and larger finite networks.
Given a finite subset $S$ of the vertices of a network,
define $\cutnetwork (S)$ to be the network gotten by throwing
away all vertices that lie outside of $S$,
along with all edges that are incident with the discarded vertices.
Define $\shortnetwork (S)$ to be the network gotten by
identifying all the vertices outside of $S$ to yield a single
new vertex $\infty$.
Call a sequence
each member of which is a finite subset $T$ of the vertices
a {\em swelling sequence}
if any finite subset $S$ of the vertices is eventually
contained in $T$.
We will be particularly interested in properties of the infinite network
that can be expressed as limits for swelling $T$ of properties of
$\cutnetwork (T)$ or $\shortnetwork (T)$.

\section{The even and odd flows.}

Given an infinite network,
and any flow through the edges of the network,
we can compute the source strength
at each of the vertices of the network,
just as we would for a finite network.
Conversely, given a source distribution,
we would like to associate a flow to it.
We will assume that there are only a finite number of
vertices where the source strength is non-zero.
For the moment, we will also assume that
the source distribution is {\em balanced},
that is,
that the total source strength is zero.

Given a balanced finite source distribution,
there are two canonical associated flows,
which we will call the {\em even flow} and the {\em odd flow}.
The even flow corresponds intuitively to a network
whose boundary at infinity is insulated.
It is the limit of the flows you get by cutting the network
down to a finite subset of the vertices,
and then letting the finite subset swell to fill up the
whole infinite network.
The odd flow corresponds to a network that is shorted together
at infinity.
It is the limit of the flows you get by shorting together
all the nodes outside of a finite subset.

The odd flow can be characterized as the unique flow of minimum energy
among all flows having the specified source distribution.
The dual characterization of the even flow
describes it in terms of functions having prescribed values,
which is not quite what we want.
Instead,
as noted by Flanders
\cite{flanders:infinite},
the even flow can be characterized
as the unique flow of minimum energy
among all flows that are limits
(in the energy dissipation norm)
of flows having the specified source distribution.

The names even and odd come from the method of images
in classical electrical theory.
Imagine doubling the network by taking two copies of it and glueing
them together along their ideal boundaries at infinity---whatever
that might mean---just as you
would double the unit disk to get a sphere.
We will call the two sides of the double
the {\em bright side} and the {\em dark side}.
The even flow corresponds to extending the source distribution
on the bright side symmetrically to the dark side,
the odd flow to extending it antisymmetrically.

\section{The even and odd resistances.}

Given any two distinct points $p$ and $q$,
place a source of strength $+1$ at $p$
and a source of strength $-1$ at $q$.
Call the resulting flows the even and odd flows from $p$ to $q$.
Call the dissipation of the even flow the even resistance
between $p$ and $q$,
and similarly for the odd resistance.

The odd flow from $p$ to $q$ is the unique unit flow from $p$ to $q$
of minimum energy,
so the odd resistance is always less than or equal to the even resistance,
and the two resistances are equal if and only if
the even and odd flows agree.
As before, the even flow is the limit of flows of compact support,
and minimizes energy among unit flows from $p$ to $q$ with this
property.
It is also characterized as the unit flow that is proportional
to the flow of the unique harmonic function that has minimum
energy among all functions taking values $1$ at p and $0$ at $q$.

\section{Shorting and cutting.}

For finite networks,
Rayleigh's cutting law states that for any two vertices $p$ and $q$,
cutting away part of the network can only make the resistance between
$p$ and $q$ bigger.
By the same token, shorting parts of the network together can only
make the resistance smaller.

These laws continue to hold for infinite networks
to the fullest extent possible.
In particular, for any finite subset $S$ containing $p$ and $q$
we have
{\samepage
\begin{eqnarray*}
\resistance (p,q,\shortnetwork (S))
\leq
\oddresistance (p,q)
\\
\leq
\evenresistance (p.q)
\leq
\resistance (p,q,\cutnetwork (S))
.
\end{eqnarray*}
}

\section{The ghost flow.}

For any balanced finite source distribution,
the difference between the even flow and the odd flow
is the flow of a harmonic function,
and in particular is sourceless.
It is the flow on the bright side of the double that results from
placing sources of twice the specified strength
at the corresponding places on the dark side.
Let's call this flow the {\em ghost flow}.

The ghost flow ``from $p$ to $q$'' has the following remarkable properties:
Its energy is the difference of the even and odd resistances,
so it vanishes if and only if the even and odd resistances agree.
Among all flows of harmonic functions $u$ of finite energy,
the ghost flow is the unique flow that minimizes the quantity
\[
\energy (u) - (u(p)-u(q))
.
\]
Thus if all finite energy harmonic functions $u$ have
$u(p)=u(q)$,
and in particular if the network has no non-constant
finite energy harmonic functions,
then the ghost flow vanishes.
Otherwise,
there is an essentially unique harmonic function of unit energy
that maximizes the difference
$u(p)-u(q)$,
and the ghost flow is proportional to the flow of this function.

From all of this,
we conclude that the even and odd resistances between $p$ and $q$ agree
if and only if
all finite energy harmonic functions $u$ have
$u(p)=u(q)$,
and that the even and odd resistances agree
for all pairs of points $p$ and $q$
if and only if
there are no non-constant finite energy harmonic functions,
Note that to check that all of the even and odd resistances agree,
we need only check the cases where $p$ and $q$ are adjacent.

\section{Networks with no non-constant finite energy harmonic functions.}

A network having no non-constant finite energy harmonic functions
is said to belong to the class $\ohd$,
or to be $\ohd$.
For an $\ohd$ network all even and odd resistances agree,
and by superposition,
the even and odd flows agree for any balanced finite source distribution.
This property characterized $\ohd$ networks,
so if you want to avoid getting into arguments about how to define
the flow corresponding to a specified source distribution,
you had best stick to $\ohd$ networks.

The class $\ohd$ fits into the hierarchy
\[
\og \subset \ohp \subset \ohb \subset \ohd
,
\]
Here $P$ and $B$ stand for `positive' and `bounded,'
and $\og$ refers to the class of surfaces
having no positive Green's function.
A network is $\og$ if and only if it has
no flow out to infinity of finite energy,
by which we mean a flow of finite energy whose source strength
is non-negative everywhere and positive somewhere.
If a network is $\og$ we call it {\em recurrent};
otherwise we call it {\em transient}.
The terminology comes from probability theory:
If you carry out a random walk on the vertices
of the network,
where at each step
you walk along one of the edges leaving the vertex you're at,
with probability proportional to the conductance of the edge,
then on a recurrent network you're certain to return eventually to your
starting point,
but on a transient network there is a positive probability that you
will wander off and never return.

The inclusions in this hierarchy imply that if you know that a network is
recurrent, or at least that it has no non-constant positive harmonic
functions,
then you know that it has no non-constant finite energy harmonic functions,
and thus that even and odd flows and resistance always agree.

\section{Flows to infinity.}

So far,
we have been dealing only with source distributions that are balanced.
For transient networks
we can relax this condition,
though now we can only consider the odd flow.
As before,
this flow is the limit of flows in the networks obtained by shorting
together nodes outside of a finite set.
As before,
it is the unique energy minimizing flow having the specified
source distribution.

\section{The $d$-grid.}

From probability theory,
we have the following standard facts
(see Spitzer
\cite{spitzer:principles},
and also Avez
\cite{avez:ohp}):
For $d=1,2$, the $d$-grid is recurrent,
i.e. $\og$.
For $d \geq 3$ the $d$-grid is transient,
but lies in $\ohp - \og$.
It follows that any grid is $\ohd$.
(We will give an independent proof of this fact later on.)
Thus for the $d$-grid 
even and odd resistances agree,
and there is no argument about what resistance means.

\section{Symmetry and superposition.}

The symmetry-and-superposition computation
of the effective resistance across
an edge of the $d$-grid depends on being able to represent the
unit flow from $p$ to $q$ as the superposition of symmetrical
flows out of $p$ and into $q$.
If $d \geq 3$,
so that the grid is transient,
we can write the flow from $p$ to $q$ as the superposition
of the odd unit flow out of $p$ and the odd unit flow into $q$,
and these flows must be symmetrical by uniqueness.
So far as the case $d \geq 3$ goes,
then,
the engineers are completely vindicated.

\section{The 2-grid.}

That leaves the $2$-grid.
Why can the flow from $p$ to $q$ be written as the superposition
of symmetrical flows out of $p$ and into $q$?
Again, the engineers' explanations vary,
but they clearly all believe that
even though the network is recurrent you can define a unique unit current
flow from any vertex $p$ out to infinity as a limit of the flows
you get as follows:
Choose a large finite subset of the vertices,
and define a source distribution
consisting of a source of strength $+1$ at $p$ and a finite
number of negative sources (i.e. sinks) at points outside of the
specified finite subset.
As the finite subset swells, the flow thus determined should
converge to a uniquely determined flow from $p$ to infinity.

Let's say that a network having this property
has a {\em good enough Green's function}.
On such a network you can work with flows out to infinity
much as you would on a transient network,
where you have a bona fide Green's function.
In particular,
you can get the flow from $p$ to infinity as the limit for swelling $T$
of the flow from $p$ to $\infty$ in $\shortnetwork (T)$,
or as the limit as $q$ marches off to infinity
of the flow from $p$ to $q$ in the infinite graph.

Note that the $1$-grid does not have a good enough Green's function:
If you short together the nodes outside a finite set,
and let the set go to infinity in a lopsided way,
you get a lopsided limit,
if indeed you get a limit at all.
Similarly, the flow from $p$ to $q$ changes abruptly when
$q$ jumps from one side of $p$ to the other.

For the $2$-grid you can actually write down the Green's function,
and then make explicit estimates to show that it is good enough.
This was the method used by Flanders
\cite{flanders:infiniteGrid}
to demonstrate
that it is valid to apply the symmetry-and-superposition method
to the $2$-grid.
This is all right as far as it goes,
but it would be nice to have a more conceptual proof,
based on the analogy between the grid and the plane.
After all, the reason everyone believes that the $2$-grid has
a good enough Green's function is that the grid looks just like the plane,
and the plane has a good enough Green's function.
If we can somehow make sense of this argument,
then without further ado we ought to be able to carry over
our results for the $2$-grid
to the symmetrical graphs shown in Figure 2.
\ldots

\fig{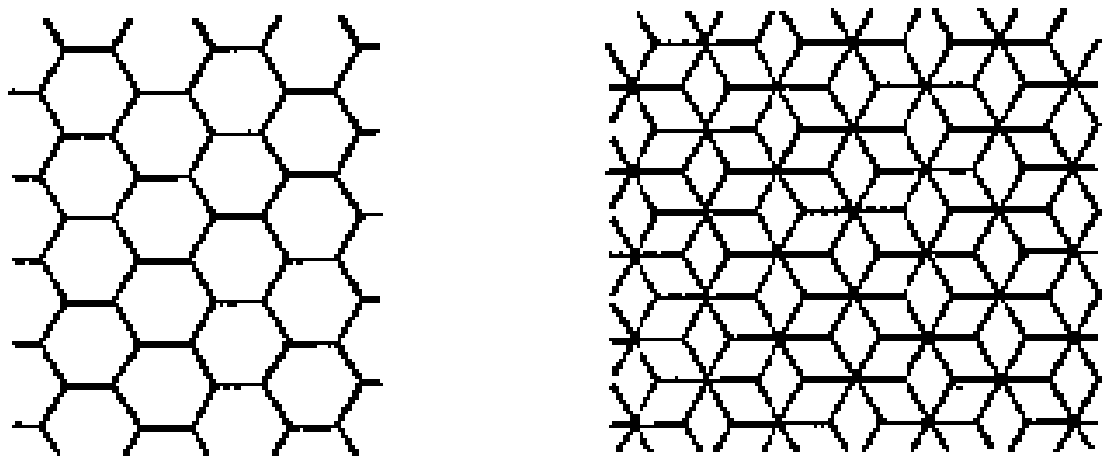}{Symmetrical graphs that look just like the plane.}

\section{Foster's method.}

We turn now to Foster's method of computing
the resistance across an edge of the $d$-grid.
Thomassen
\cite{thomassen:infinite}
shows how to make rigorous the application of
Foster's theorem for finite networks to smallish symmetrical infinite
networks like the $d$-grid.
We will review Thomassen's analysis,
and extend his results in the light of the network theory
we have developed.

Given a finite subset $S$ of the vertices of a graph $G$,
define $\boundary (S)$ to consist of those vertices of $S$ that
are connected by an edge to a vertex outside $S$.
Define $\edgeboundary (S)$
to consist of all edges connecting a vertex in $S$ to a vertex
outside $S$.
Call a swelling sequence along which
\[
{{|\edgeboundary (T)|}\over{|T|}}
\]
goes to $0$ a {\em pinching sequence}.
Call a graph {\em smallish} if it has
a pinching sequence.
We will only be considering graphs of bounded valence;
such a graph is smallish
if and only if it has a swelling sequence
along which
\[
{{|\boundary (T)|}\over{|T|}}
\]
goes to $0$.

Say that two vertices $v, v'$ are {\em of the same kind}
if $v$ can be mapped to $v'$ by a symmetry of the graph.
Similarly for edges.
Call a graph {\em symmetrical}
if it has only a finite number of different kinds of vertices
(or edges---it makes no difference).

In a smallish symmetrical graph,
the vertices have well-defined relative frequencies,
which are positive rational numbers summing to $1$.
Similarly for edges.
These relative frequencies are the limits of the actual
relative frequencies
in $\cutgraph (T)$ or $\shortgraph (T)$,
where $T$ runs out along any pinching sequence.
Define the average of any function that depends only on the
kind of a vertex (or edge)
with respect to these relative frequencies.
This average will be the limit of the average in
$\cutgraph (T)$ or $\shortgraph (T)$ as $T$ runs out along
any pinching sequence.

{\bf Theorem.}
In a smallish symmetrical graph,
\begin{eqnarray*}
\averageevenresistance (G)
=
\averageoddresistance (G)
\\
=
{2 \over {\averagevalence (G)}}
.
\end{eqnarray*}

{\bf Proof.}
For any finite set of vertices $T$,
Rayleigh's cutting law implies that
the even resistance in the infinite graph
across any edge that remains in $\cutgraph (T)$ is 
less than or equal to the (plain old) resistance
across that edge in $\cutgraph (T)$.
But by Foster's theorem,
when $|T|$ is large the average resistance
across the edges of $\cutgraph (T)$
is very nearly
\[
{2 \over {\averagevalence ( \cutgraph (T) )}}
.
\]
Taking limits along a pinching sequence,
we find that
\[
\averageevenresistance (G) \leq 2 / \averagevalence (G)
.
\]
An analogous shorting argument yields
\[
\averageoddresistance (G) \geq 2 / \averagevalence (G)
.
\]
But even resistances are always at least as big as odd resistances,
so
\[
\averageevenresistance (G) \geq \averageoddresistance (G)
.
\]
This closes up the circle of inequalities,
and the theorem follows.
$\qed$

{\bf Corollary.}
In a smallish symmetrical graph $G$,
all even and odd resistances agree,
and $G$ is $\ohd$.
$\qed$

{\bf Note.}
Actually,
we expect that a smallish symmetrical graph is $\ohp$.
But Geoff Mess tells me that this isn't true
\ldots

{\bf Corollary.}
In a smallish vertex-transitive graph $G$,
the average of the resistances across the edges emanating from any given vertex
is
\[
2 / \valence (G)
. \dispqed
\]

In an edge-transitive graph, there can be either one type of vertex
or two types.
If there are two types,
their relative frequencies are inversely proportional to their valences,
which we denote by $\valenceone (G)$ and $\valencetwo (G)$.

{\bf Corollary.}
Let $G$ be a smallish edge-transitive graph.
If $G$ has one type of vertex,
the resistance across any edge is
\[
{2 \over {\valence (G)}}
.
\]
If $G$ has two types of vertices,
the resistance across any edge is
\[
{{\valenceone (G) + \valencetwo (G)} \over {\valenceone (G) \valencetwo (G)}}
. \dispqed
\]

{\bf Note.}
This result follows from the symmetry-and-superposition argument
as well,
though we still haven't justified the use of this method if the 
graph is recurrent.
We expect that a smallish symmetrical graph can only be recurrent
if it looks like the $2$-grid.
Geoff Mess tells me that this follows from the work of 
Gromov
\cite{gromov:polynomial}
and Varopoulos
\cite{varopoulos:polya}.
\ldots

\newpage

\bibliography{net}
\bibliographystyle{plain}

\end{document}